\documentclass[preprint,12pt]{article}

\usepackage{amsmath,amssymb,amsthm}
\usepackage{mathtools}
\usepackage{microtype}
\usepackage{hyperref}
\usepackage{graphicx}
\usepackage{booktabs}
\usepackage{enumitem}
\usepackage{cleveref}
\usepackage{url}
\usepackage{accents}
\usepackage{setspace}
\usepackage[a4paper, left=3.5cm, right=3.5cm]{geometry}
\usepackage{hypernat}  
\hypersetup{hypertexnames=false}
\DeclareFontFamily{OMX}{yhex}{}
\DeclareFontShape{OMX}{yhex}{m}{n}{<->yhcmex10}{}
\DeclareSymbolFont{yhlargesymbols}{OMX}{yhex}{m}{n}
\DeclareMathAccent{\wideparen}{\mathord}{yhlargesymbols}{"F3}

\newtheorem{theorem}{Theorem}
\newtheorem{proposition}[theorem]{Proposition}

\theoremstyle{definition}
\newtheorem{definition}[theorem]{Definition}

\newtheorem*{corollary*}{Corollary}

\theoremstyle{definition}
\newtheorem*{definition*}{Definition}
\newtheorem*{remark*}{Remark}

\providecommand{\keywords}[1]{\textbf{\textit{Keywords:}} #1}
\providecommand{\MSC}[1]{\textbf{\textit{MSC subject classifications:}} #1}

\begin{document}

\title{From scalar rigidity to  abstract prime set and abstract integer system}

\author{W. Oukil\thanks{Corresponding author. Email: \texttt{oukil.walid@gmail.com}} \\
\small Faculty of Mathematics. \\
\small University of Science and Technology Houari Boumediene. \\
\small BP 32 EL ALIA 16111 Bab Ezzouar, Algiers, Algeria.}

\date{\today}

\maketitle


\begin{abstract}
We propose a unified framework for Prime Rigidity Theory (PR) by introducing the notions of an "abstract integer system" and an "abstract prime set", consisting, respectively, of two unbounded sequences $\mathcal{P}$ and $\mathcal{N}$ of elements of $\mathcal{X}$, where $\mathcal{X}$ denotes the algebra of bounded linear endomorphisms of  a complex Banach space $X$, such that every element of $\mathcal{N}$ can be expressed as a finite product of powers of elements of $\mathcal{P}$. We define the zeta function $\zeta_{X,\mathcal{N}}$ on $\mathcal{X}\setminus{I_X}$, where  $I_X$ is the identity endomorphism. The scalar restriction $\zeta_{\mathbb{C},\mathcal{N}}$ yields generalized Beurling zeta functions; in particular, $\zeta_{\mathbb{C},\mathbb{N}}$ coincides with the classical Riemann zeta function. We introduce the notions of a "rigid abstract integer system" and  "rigid abstract prime set" in order to preserve the rigid asymmetric structure obtained in the scalar case.

\end{abstract}


\MSC{34A30, 11M06, 47A60, 20M05.}


\section{Introduction} 
 The asymmetry of a parametrized non-homogeneous linear complex differential equation on $[1,\infty)$, with  non-homogeneous term $f\in L^\infty([1,+\infty),\mathbb{R}_+)$ satisfying the Rotation number hypothesis, is characterized by the fact that two symmetric solutions corresponding to the parameters $s$ and $1-s$, and satisfying the same initial condition, cannot both be bounded on $[1,+\infty)$.

The space $\mathcal{X}$ denotes the  algebra of bounded linear endomorphisms of a complex Banach space $X$.  When $X=\mathbb{C}$, we can identify $End(\mathbb{C})$ with $\mathbb{C}$, where each endomorphism acts by multiplication on $\mathbb{C}$.   

The Euler product and the functional equation  can be derived by the same mechanism as in the scalar case.

We identify the Riemann zeta function $\zeta(s)$, for $s\in\mathbb{C}\setminus{1}$, with a parametrized non-homogeneous linear system with parameter $s$ and the fractional part function as its non-homogeneous term. We then replace the complex parameter $s$ by an endomorphism and the fractional part function by a souhaitabel real-valued   function $p$ depending of the "abstract prime" set $\mathcal{N}$ which is an unbounded  sequence of $\mathcal{X}$. The Euler product and the functional equation  can be derived by the same mechanism as in the scalar case. This allows us to introduce a generalization of the zeta function $\zeta_{X, \mathcal{N}}$ on the Banach space $\mathcal{X}$. 
 his leads to the following diagram:

$$
\begin{array}{ccc}
V\in  \mathcal{X}
& \xrightarrow{\quad\zeta_{X, \mathcal{N}}\quad} &
 \zeta_{X, \mathcal{N}}(V)
\\[2mm]
\downarrow\scriptstyle{  \mathcal{X}= \mathbb{C},\  \mathcal{P}:=\mathcal{N}\subset \mathbb{R} }
&&
\downarrow
\\[2mm]
s\in\mathbb C,\ \mathcal{P}
& \xrightarrow{\quad\zeta_\mathcal{P}\quad} &
\zeta_{\mathcal{P}}(s)\\ 
\\[2mm]
\downarrow\scriptstyle{  \mathcal{P} =\mathbb{N} }
&&
\downarrow
\\[2mm]
s\in\mathbb C,\ \mathbb{N}
& \xrightarrow{\quad\zeta\quad} &
\zeta(s).
\end{array}
$$
where $\zeta_\mathcal{P}:\mathbb{C}\to \mathbb{C}$ is the generalized Beurling zeta function with   prime set $\mathcal{P}\subset \mathbb{R}$.

\section{Scalar Rigidity}\label{sec:scalar}
Consider the following  non-homogeneous linear complex differential equations:
\begin{equation}\label{EDOSimplie} 
  \dot{\phi} = w \,  t^{-1}\,  \phi +(1 -w) \,  t^{-1}\,   \eta(t), \quad\phi(1)=1,\quad \phi: [1,+\infty)\to \mathbb{C},
\end{equation}
where $w \in \mathbb{C}$   is the parameter and the non-homogeneous term  $\eta:[1,+\infty)\to\mathbb{R}_+$  belongs to   $L^{\infty}([1,+\infty), \mathbb{R}_+)$. 
The unique   continuous solution $\psi_{\eta,w}:\mathbb{R}_+\to\mathbb{C}$ of the previous differential equation  is given by
\begin{equation}\label{Solutionz}
\psi_{\eta,w}( t) = t^{w}\Big[1+  (1-w)\int_0^{t} u^{-1-w}\eta(u) \, du\Big], \quad \forall t \ge 1.
\end{equation}

We  assume that the function $\eta$ satisfies the following {\it {Rotation number hypothesis}} (see \cite{Oukil}):
\begin{equation*}\label{NombreRotation}\tag{\textbf{H}}
 \exists \rho_\eta>0:\quad   \sup_{t\geq1}\Big| \int_1^t\, \big( \eta(u)-\rho_\eta\big)\, du\Big|<+\infty,
\end{equation*}   
\vspace{0.1cm}
The  number $\rho_\eta$  is called the {\it rotation number of}  $\eta$.
We denote by $B \subset \mathbb{C}$,   the critical strip except the critical line, as
\begin{equation}\notag
B := \Big\{ w \in \mathbb{C} : \quad\Re(w)\in(  0,1)\setminus\{\frac{1}{2}\} \Big\}.
\end{equation}  
Next, we consider the following statement, stated in\cite{Oukil}, as a characterization rather than as a result.

\begin{theorem}\label{MainTheo}
Let $\eta \in L^\infty([1,+\infty),\mathbb{R}_+)$ satisfy \eqref{NombreRotation}. Then for every $s\in B$ we have the following implication:
\[
\sup_{t\geq0} \| \psi_{\eta,s}(t) \|<+\infty   \implies\sup_{t\geq0}\| \psi_{\eta,1-s}(t) \|=+\infty,
\]
where $\psi_{\eta,w}$ is the unique   continuous solution of the differential equation \eqref{EDOSimplie}.
\end{theorem}
The fractional part function $\eta_*(t) = \{t\}$ is bounded and locally integrable. Since the function $p(t) = \int_1^t (1/2 - \{u\})\,du$ is $1$-periodic and vanishes at integer values, it is bounded for all $t \ge 1$. Consequently, the function $\eta_*(t) = \{t\}$ satisfies hypothesis \eqref{NombreRotation} with $\rho = \frac{1}{2}$. More precisely, every nonzero periodic function in $L^\infty([1,+\infty),\mathbb{R}_+)$  satisfies  the hypothesis \eqref{NombreRotation}. From the integral representation of the Riemann zeta function (Titchmarsh, \cite{Titchmarsh}, page 14, Equation 2.1.5):
\begin{equation}\label{Zetafunction}
\forall w \in \mathbb{C}_+ \setminus \{1\}: \quad \frac{1-w}{w}\zeta(w) =  \mu_{\eta_*}(w).
\end{equation}
\section{Abstract integer system and abstract prime set}
 When we restrict the Banach space to the complex scalar field $\mathcal{X} = \mathbb{C}$, the following "abstract prime set" are reduced to the Beurling's generalized prime\cite{DiamondZhang2016}.
\begin{definition}
 An \emph{abstract integer system} is given by the couple $(X, \mathcal{N})$ where $X$ is a complex Banach space and $\mathcal{N}$ is  an unbounded and strictly increasing sequence of element of $\mathcal{X}$  such that there exists a sequence $(x_n)_{n\in \mathbb{N}^*} \subset \mathcal{X}$ with
$$\mathcal{N}:=  \Big\{\prod_{(m,k)\in A} x_m^{k}:\quad A\subset \mathbb{N}^*\times \mathbb{N} \Big\}$$.\\
We call the unbounded sequence $\mathcal{P}:=(x_n)_{n\in \mathbb{N}^*}$ the  {\it{abstract prime set}} of the abstract integer system  $(X, \mathcal{N})$.
\end{definition}
\subsection{Zeta function over endomorphisms set} 
Let $X$ be a complex Banach space.  We denote by $\mathcal{X}$ the space of all bounded linear endomorphisms from $X$,  into itself. For $T\in \mathcal{X}$, the operator norm is defined by
\[
\|T\|:=\sup_{\|x\|\le 1}\|T\, x\|.
\]   
With this norm, $\mathcal{X}$ is itself a complex Banach space. 
 For $W \in \mathcal{X}$, we denote by $\sigma(W)$ the spectrum of the endomorphism $W$. We define the real number  $\alpha(W)$ as 
\begin{equation}\label{DefalphaW}
\alpha(W):=\inf_{\lambda\in \sigma(W)}\Re(\lambda).
\end{equation}
Let $( X, \mathcal{N})$ be an abstract integer system  and define
\[
\mathcal{B}_{\alpha_{\mathcal{N}}}:=\big\{W\in \mathcal{X}:\quad \alpha(W)>\alpha_{\mathcal{N}}\big\}.
\]
where  the positive  constant $\alpha_{\mathcal{N}}>0$ is defined as
\[
\alpha_{\mathcal{N}}=\inf\Big\{\alpha(W) :\quad \sum_{v\in \mathcal{N}}\frac{1}{\|v\|^{\alpha(W)}}<+\infty,\quad W\in \mathcal{X}\Big\}.
\]
We define the generalized zeta function $\zeta_{X,\mathcal{N}}(V): \mathcal{X}\to \mathcal{X}$ as following
\[
\zeta_{X,\mathcal{N}}(V):=\sum_{v\in \mathcal{N}}\frac{1}{\|v\|^V},\quad V\in \mathcal{B}_{\alpha_{\mathcal{N}}}.
\] 
Use the following commutativity $\|x\|\, \|y\|=\|y\| \, \|x\|$, rhe   Euler product is given by 
\[
\sum_{v\in \mathcal{N}}\frac{1}{\|v\|^V}= \prod_{v\in \mathcal{P}} \, \frac{1}{ 1-\|p\|^{-V}},\quad V\in \mathcal{B}_{\alpha_{\mathcal{N}}},
\] 
where $\mathcal{P}$ is the abstract prime set of $(X, \mathcal{N})$.
\section{Rigid abstract  integer system}
Define 
\[
\Psi[g](t):=\int_1^t \, u^{-1}\, g(u)\, du,\quad \forall  g\in L_{loc}^\infty([1,+\infty),\mathbb{R}_+).
\]
Let   $\mathcal{N}:=(v_n)_{n\in \mathbb{N}}$ be an  abstract integer system of  $X$. Denote
\begin{equation}\label{Cunting}
N(t)=\sum_{\|v_n\|\le t} 1
\end{equation}
We call the previous function the {\it{ the counting function }} of $(X, \mathcal{N})$. As in scalar case, let $c\in\mathbb{R}$ and   $W\in \mathcal{X}$ define  ${\Phi}_W:[1,+\infty)\to \mathcal{X}$ as
\begin{equation}\label{ExpliEDO}
{\Phi}_{c,N}[W](t)=c\, t^{W}+\sum_{k=0}^{+\infty}\,\Psi^{k} [p ](t)  \big(W^{k}-W^{k+1 }\big);\  p(t):=c\, t- N(t),\ \forall t\geq1,
\end{equation} 
where
\[
t^{W}:= \exp\big(W\, \ln(t)\big)=\sum_{n=1}^{+\infty}\frac{1}{n!}\, \ln^n(t)\, W^n.
\]
The previous function is the unique continuous solution of the following linear differential equation with non-homogeneous term $\eta$:
\begin{equation}\label{endomorphismEDO} 
\frac{d}{dt}x(t)=W \, t^{-1}\, x(t)+(I-W) \, t^{-1}\, \big(c\, t- N(t)\big),\quad t\geq1,\quad x(1)= c\, I_X .
\end{equation}
where $I_X$ is the identity endomorphism of $X$.    We denote critical endomorphisms strip $\mathcal{B}(X)\subset \mathcal{X}$ as
\[
\mathcal{B}(X):=\big\{W\in \mathcal{X}:\quad \alpha(W)\in(0,1)\big\}
\]
where we recall that    $\alpha(W)$ is given by\eqref{DefalphaW}.
\begin{definition}\label{defEndfSpace}
Le $c\in (0,1]$ and let $(X, \mathcal{N})$ be an abstract integer system with counting function $N$ as defined \eqref{Cunting}.     We say that $(X, \mathcal{N})$  is a \emph{rigid  abstract integer system with slope} $c$ if
\begin{itemize}
\item   the real function $\eta(t):= ct-N(t)$ for all   belong to $L^{\infty}([1,+\infty), \mathbb{R})$;
\item   for every $W\in \mathcal{B}(X)$  we have the following implication
$$
\sup_{t\geq0}\big\|{\Phi}_{c, N}[W](t)\big\|<+\infty \implies \sup_{t\geq0}\big\|{\Phi}_{c, N}[I-W](t)\big\|=+\infty.
$$
\end{itemize} 
We  call the abstract prime set of  $(X, \mathcal{N})$ the   \emph{rigid abstract prime}. 
\end{definition}
The primes  set $\{2,\ 3,\ 5,\, ...\}$ of $\mathbb{N}$ is contained on the set of  the "rigid abstract prime set"  when    we restrict the Banach space to the complex scalar field $X = \mathbb{C}$.

 \section{Zeta function over   $End^+(X)$}

For an endomorphism $V\in \mathcal{X}$ of a Banach space, we have
 \[
u^{-1-V}:=u^{-1}\, \exp\big(-V\, \ln(u)\big)=\sum_{n=1}^{+\infty}\frac{(-1)^n}{n!}\, \ln^n(u)\, V^n.
\]
Denote
\[
End^+(X):=\big\{W\in \mathcal{X}\setminus\{I\}:\quad  \alpha(W)>0\big\}.
\] 
where $\alpha(W)$ is as given in \eqref{DefalphaW}.
\begin{proposition}\label{DefZetap}[Abel representation of $\zeta_{X,\mathcal{N}}$]
Let $( X, \mathcal{N})$ be a rigid abstract integer system  with counting function $N$ and slope $c\in \mathbb{R}$. 
\[
\zeta_{X,\mathcal{N}}(V)= -V  \int_{0}^{+\infty} u^{-1-V} \, \big(c\, u-N(u)\big) \, du,\quad \forall V\in End^+(X).
\]
\end{proposition}

\subsection{Functional equation} 
Using the change of variable $\|v \| t\to \mu$, we give in the following theorem the  functional equation
\begin{proposition}
Let $( X, \mathcal{N})$ be a rigid abstract integer system with counting function $N$ and slope $c\in \mathbb{R}$. Suppose that there exists a function $h:\mathbb{R}\to \mathbb{R}$ such that
$$  N(t)=c\, t-\sum _{v\in \mathcal{N}}h\big( \|v \| t\big)$$
Then
\[
\zeta_{X, \mathcal{N}}(V)=\zeta_{X, \mathcal{N}}(I-V)\int_0^{+\infty} u^{-1-V} \, h(u)\, du,\quad \forall V\in End^+(X).
\]   
\end{proposition}
\section{Conclusion}

The aim was to construct a framework in which all propositions and results are established in the same way as in the scalar case, based on the underlying differential equation. The approximations of $\pi(x)$ are also preserved.

\end{document}